\documentclass{article}
\usepackage{a4wide} 
\usepackage{amssymb,amsmath}
\usepackage{graphicx}

\newtheorem{proposition}{Proposition}
\newtheorem{lemma}{Lemma}
\newtheorem{definition}{Definition}

\newcommand{\ts}{\hspace{0.5pt}}
\newcommand{\nts}{\hspace{-0.5pt}}

\newcommand{\ZZ}{\mathbb{Z}}
\newcommand{\RR}{\mathbb{R}\ts}
\newcommand{\CC}{\mathbb{C}\ts}
\newcommand{\NN}{\mathbb{N}}
\newcommand{\XX}{\mathbb{X}}
\newcommand{\vG}{\varGamma}
\newcommand{\vL}{\varLambda}
\newcommand{\per}{\mathrm{per}\ts}

\begin{document}

\begin{center}

\begin{Large}\begin{bf}
On the Notions of Symmetry and Aperiodicity\\[1ex] for Delone Sets\vspace{3ex}
\end{bf}\end{Large}

\begin{large}
Michael Baake$^{1}$ and Uwe Grimm$^{2}$\vspace{3ex}
\end{large}

\begin{small}
$^{1}$Fakult\"{a}t f\"{u}r Mathematik, Universit\"{a}t Bielefeld,\\
       Postfach 100131, 33501 Bielefeld, Germany\\[0.5ex]
$^{2}$Department of Mathematics and Statistics, The Open University,\\
       Walton Hall, Milton Keynes MK7 6AA, United Kingdom\vspace{1ex}
\end{small}
\end{center}

\begin{quote}Non-periodic systems have become more important in recent
  years, both theoretically and practically. Their description via
  Delone sets requires the extension of many standard concepts of
  crystallography. Here, we summarise some useful notions of symmetry
  and aperiodicity, with special focus on the concept of the hull of a
  Delone set. Our aim is to contribute to a more systematic and
  consistent use of the different notions.
\end{quote}

\section{Introduction}

The notion of symmetry for (fully) periodic objects such as perfect
crystals is well established and has a long history.  While there are
slightly different approaches, for instance via direct or via Fourier
space, there is little dispute what the symmetry of a given structure
really is.  As it turns out, an important ingredient is the underlying
periodicity, which brings in a deep connection with lattices and their
structures; see \cite{Schw} for a systematic mathematical exposition.

Things are more complex in the setting of general (discrete)
structures, in particular in the absence of non-trivial periods. In
particular, rather different notions are in use of what terms such as
`periodicity' or `symmetry' are supposed to mean. While this should be
no problem in principle, experience proves otherwise. Quite some
confusion emerges from the comparison of `results' that were phrased
in different contexts, which asks for some clarification, as has been
noted in this context in \cite[Ch.~10]{GS}.

The purpose of this article is to contribute to this question, via
some notions that have proved particularly useful in the setting of
non-periodic and aperiodic structures; see \cite{Pat,BMBook} for
background. We present some definitions and results, all
example-driven and designed for Euclidean space $\RR^{d}$, and add
some arguments in favour of their systematic use, one being a higher
degree of consistency.

In view of the well-established equivalence concept of mutual local
derivability (MLD), compare \cite{B,BSJ} and the discussion below, we
do this mainly for Delone sets, and leave it to the reader to rephrase
results in terms of his or her favourite structure, be it a tiling, a
marked pattern, an ornament or any other MLD representative. Some of
our examples will also be drawn from the tiling zoo, but we primarily
view them as Delone sets (for instance, via their vertex points) in
the spirit of \cite{Lag99}.

Very little (if anything) of what we discuss below is really
new. However, the bits and pieces are scattered over the literature,
both in mathematics and physics. Due to the renewed interest in
non-periodic systems following the award of the 2011 Nobel Prize in
Chemistry to Dan Shechtman for his discovery of quasicrystals
\cite{Danny}, it seems necessary to recapitulate and slightly adjust
some of the most relevant concepts in a more coherent fashion. We hope
that our brief exposition below will help to clarify the notions and
avoid future misunderstandings.

\section{General notions and definitions}

A central concept for our discussion is that of a \emph{Delone set}
$\vL \subset \RR^{d}$, by which we mean a point set (a countable
subset of $\RR^{d}$) that is neither `too crowded' nor `too sparse'.
More precisely, there are positive radii $r$ and $R$ such that no open
ball $B_{r} (x)$ contains more than one point of $\vL$, irrespective
of the centre $x$, while each ball $B_{R} (x)$ contains at least one
point of $\vL$.  More generally, one can consider other point sets as
well. In view of the applications in crystallography, we often
consider point sets that are at least \emph{uniformly discrete}, which
is the above condition with the ball $B_{r} (x)$ alone, while we might
sometimes omit the \emph{relative denseness} (which is the other
condition). A reasonable minimal assumption is the \emph{local
  finiteness} of $\vL$, by which we mean that $\vL\cap K$ is a finite
set for any compact $K\subset\RR^{d}$. Local finiteness of $\vL$ is
equivalent to $\vL$ being discrete (each $x\in\vL$ possesses an open
neighbourhood that contains no other point of $\vL$) and closed. Note
that the set $\{\frac{1}{n}\mid n\in\NN\}$ is discrete, but not closed
(and clearly not locally finite).

A Delone (or, more generally, a locally finite) set $\vL$ is said to
have a non-trivial \emph{period} when $t+\vL = \vL$ holds for some $0\ne
t\in\RR^{d}$, and \emph{non-periodic} otherwise, where
$t+\vL=\{t+x\mid x\in\vL\}$. The set
\[
    \per (\vL) \, = \, \{ t\in \RR^{d} \mid t+\vL = \vL \}
\]
is the set of \emph{periods} of $\vL$, which is always a subgroup of
$\RR^{d}$. Clearly, $\vL$ is non-periodic if and only if $\per (\vL) =
\{ 0 \}$, which means that $\per(\vL)$ is the trivial group. When $t$
is a non-trivial period of $\vL$, $\per (\vL)$ contains the subgroup
$t\ts\ZZ$. A Delone set $\vL$ is called \emph{crystallographic} (or fully
periodic) when the $\RR$-span of $\per (\vL)$ is $\RR^{d}$. The
simplest examples are \emph{lattices}, by which we mean discrete
co-compact subgroups of $\RR^{d}$. They are thus point sets of the
form
\[
    \vG \, = \, \ZZ b_{1} \oplus \ldots \oplus \ZZ b_{d}
\]
where $\{ b_{1}, \ldots , b_{d}\}$ is a basis of $\RR^{d}$.

Let us emphasise that, as soon as $d>1$, a Delone set can be periodic
without being crystallographic (some authors call this sub-periodic,
though we will not use this term). If $\vL\in\RR^{d}$ is locally
finite and $\per(\vL)$ is a (discrete) group of rank $r<d$, $\vL$ is
said to be \emph{periodic of rank $r$}. In particular,
non-crystallographic (for $d\ge 2$) is weaker than
non-periodic. Later, we shall also introduce the term
\emph{aperiodic}, which refers to a yet stronger property than
non-periodicity. Our notion of aperiodicity will be consistent
with the one introduced in \cite{GS}.

The symmetries of crystallographic point sets are fully understood,
and completely classified in small dimensions \cite{Schw,Brown}.
Particularly simple are lattices, and it may be this widely discussed
class that makes the adjustment of suitable symmetry concepts to other
(and in particular non-periodic) cases so unfamiliar. Let us thus
first look at symmetries under (linear) isometries in more detail.

\section{Rotations and reflections of individual Delone sets}

In one dimension, the reflection $r_{x} \! : \, \RR
\longrightarrow \RR$ in the point $x$ is defined by $r_{x} (y) =
2x-y$. Clearly, $(r_{x})^{2}$ is the identity on $\RR$.  A point set
$\vL\subset\RR$ is called \emph{reflection
  symmetric} in the point $x$, if
\[
      r_{x} (\vL) \, = \, \vL \ts .
\]
Observing $(r_{x} \circ r_{y} ) (z) = 2 (x-y) + z$, which means that
the product of two reflections in $\RR$ is a translation, the
existence of more than one reflection centre has the following
consequence.
\begin{lemma}\label{lem:cox1d}
  Let\/ $\vL\subset\RR$ be a point set that is reflection symmetric in
  the distinct points\/ $x$ and\/ $y$. Then, $\vL$ is periodic with
  period\/ $\ts 2\ts \lvert x-y \rvert$. \hfill $\square$
\end{lemma}

If $\vL\subset\RR$ is locally finite, the existence of two distinct
reflection centres implies, via Lemma~\ref{lem:cox1d}, the existence
of a non-trivial period, wherefore $\vL$ must also be uniformly
discrete. Consequently, there are also two reflection centres of
minimal distance, $x'$ and $y'$ say. In this case, $\vL$ is
crystallographic, with $\per (\vL) = 2 \ts \lvert x' \! - \nts y'
\rvert \ts \ZZ$ as its lattice of periods.  Moreover, the set $\bigl(x' +
\per (\vL)\bigr)\cup\bigl(y'+ \per (\vL)\bigr)$ comprises all
reflection centres of $\vL$.

In dimension $d$, one inherits this structure when the point set $\vL$
is reflection symmetric in two parallel, but distinct, co-dimension
one hyperplanes; see \cite{HumBook} for background material, including
a classification of finite and affine reflection groups.

\begin{figure}[t]
\centerline{\includegraphics[width=0.45\textwidth]{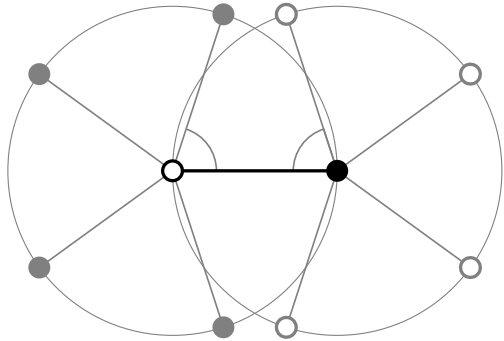}}
\caption{Sketch of the incompatibility between points of fivefold
  rotational symmetry at minimal distance (black line) and 
  periodicity in the Euclidean plane.}
\label{fig:nofivegen}
\end{figure}

Let us now take a closer look at individual rotation symmetries in the
Euclidean plane. Our (actually rather classic) arguments will revolve
around the consequences of multiple symmetry centres. For convenience,
we observe $\RR^{2}\simeq \CC$ and use complex numbers. Consider a
uniformly discrete point set $\vL\subset\CC$ and assume that it has
two distinct fivefold rotation centres, $z$ and $z'$ say, which means
that
\[
    \xi^{}_{5} (\vL - z) = \vL - z
    \quad \text{and} \quad
    \xi^{}_{5} (\vL - z') = \vL - z' ,
\]
where $\xi^{}_{5}$ is a primitive fifth root of unity. Note that
either rotation centre may be, but need not be, an element of
$\vL$. Still, the uniform discreteness of $\vL$ implies that there is
such a pair $z\ne z'$ with $\lvert z'\nts - z\rvert$
minimal. Figure~\ref{fig:nofivegen} shows how each rotation centre
produces five centres of the other type, namely $\{ z + \xi^{\ell}_{5}
(z'-z) \}$ and $\{ z' + \xi^{\ell}_{5} (z-z') \}$ with $0\le \ell \le
4$. Here, pairs of opposite centres have distance
\[
   \big\lvert \xi^{\ell}_{5} + \xi^{-\ell}_{5} - 1 \big\rvert
   \, \lvert z' -z \rvert \ts ,
\]
where the first factor takes a non-zero value $<1$ for some
$\ell$. This contradicts the separation assumption on the position of
the rotation centres, wherefore we conclude that two rotation centres
are impossible. This also implies that a uniformly discrete point set
$\vL \subset \mathbb{R}^{2}$ with fivefold symmetry cannot possess
a non-trivial period. This is the reason why fivefold symmetry is
called a \emph{non-crystallographic} symmetry of the plane; compare
\cite{CMP98,HumBook}.

The above argument applies to all other non-crystallographic
rotation symmetries ($n=5$ or $n\ge 7$) of the plane in an analogous
way. In contrast, the crystallographic symmetries `escape' this
argument, see also \cite[Sec.~6.6]{HumBook}, because the prefactor is
either always $\ge 1$ (for $1 \le n \le 4$) or the only value $<1$ is
zero (which happens for $n=6$). Together, this gives the following
result.

\begin{proposition}\label{prop:symmcent}
  Let\/ $\vL\subset\CC\simeq \RR^{2}$ be a uniformly discrete point
  set with an exact\/ $n$-fold rotational symmetry. If\/ $n$
  is non-crystallographic, which means\/ $n=5$ or\/ $n\ge 7$,
  there can only be one such rotation centre. When\/ $n\in \{3,4,6\}$,
  the existence of more than one rotation centre is possible,
  and then implies lattice periodicity of\/ $\vL$. When\/ $n=2$,
  the existence of another rotation centre means that\/ $\vL$
  is at least rank-$1$\/ periodic.
\end{proposition}
\noindent \textsc{Proof.}  The claim about the impossibility of more
than one rotation centre for non-crystallo\-graphic values of $n$
follows from our above arguments. The case $n=4$ (in standard
orientation) leads to rotation centres on the lattice $\vG = \lvert
z'\nts -z\rvert \, \ZZ^{2}$, which results in $2 \vG \subset \per
(\vL)$. Similarly, $n=6$ produces a triangular lattice, while the
choice $n=3$ generates rotation centres on the vertices of a honeycomb
packing. Finally, $n=2$ produces rotation centres along a line.  This
is (via an embedding) a consequence of
Lemma~\ref{lem:cox1d}.\hfill$\square$\bigskip

  For a systematic discussion of non-crystallographic root systems and
  their reflections in higher dimensions, we refer to \cite{CMP98} and
  references therein in conjunction with \cite{HumBook}. Let us just
  mention that a uniformly discrete point set in $3$-space can at most
  have one symmetry centre for icosahedral symmetry, which is the
  most relevant case for quasicrystals in $3$-space.
\medskip

\begin{figure}[t]
\centerline{\includegraphics[width=0.9\textwidth]{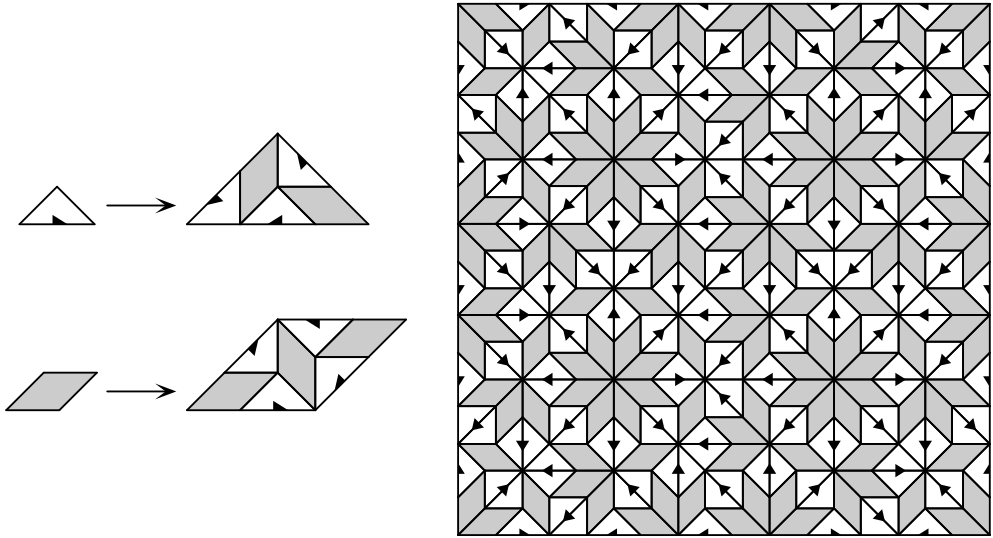}}
\caption{Inflation rule for the (essentially undecorated)
  Ammann-Beenker tiling (left), and a patch (right) obtained by three
  inflation steps from a square-shaped patch (consisting of two
  triangles) in the centre; see text for details. The patch has no
  rotation symmetry, but is reflection symmetric in the diagonal.}
\label{fig:abinfltil}
\end{figure}

  The exact symmetry of an individual point set or pattern is an
  important concept, but it is of limited use for non-periodic
  structures. Let us briefly explain this well-known observation with
  a classic example. We employ an inflation rule for a tiling made
  from a square (which is dissected into two triangles) and a rhombus,
  as displayed in Figure~\ref{fig:abinfltil}. It is known as the
  (undecorated) Ammann-Beenker tiling; see \cite{GS,AGS,TAO} and
  references therein for details. In an inflation step, a tiling is
  scaled by a factor (in this case $1+\sqrt{2}$) and the inflated
  tiles are then dissected into copies of the original size, as
  indicated in the left panel of Figure~\ref{fig:abinfltil}. Iterating
  this procedure on an appropriate initial patch (which may consist of
  a single tile) will lead to a space-filling tiling of the plane.  

    The patch shown in Figure~\ref{fig:abinfltil} is part of a
  tiling of the plane that is a fixed point under the square of the
  inflation rule (with a single square as its central seed).  It has
  no rotational symmetry at all, and a reflection symmetry only in the
  diagonal.

\begin{figure}[t]
\centerline{\includegraphics[width=0.8\textwidth]{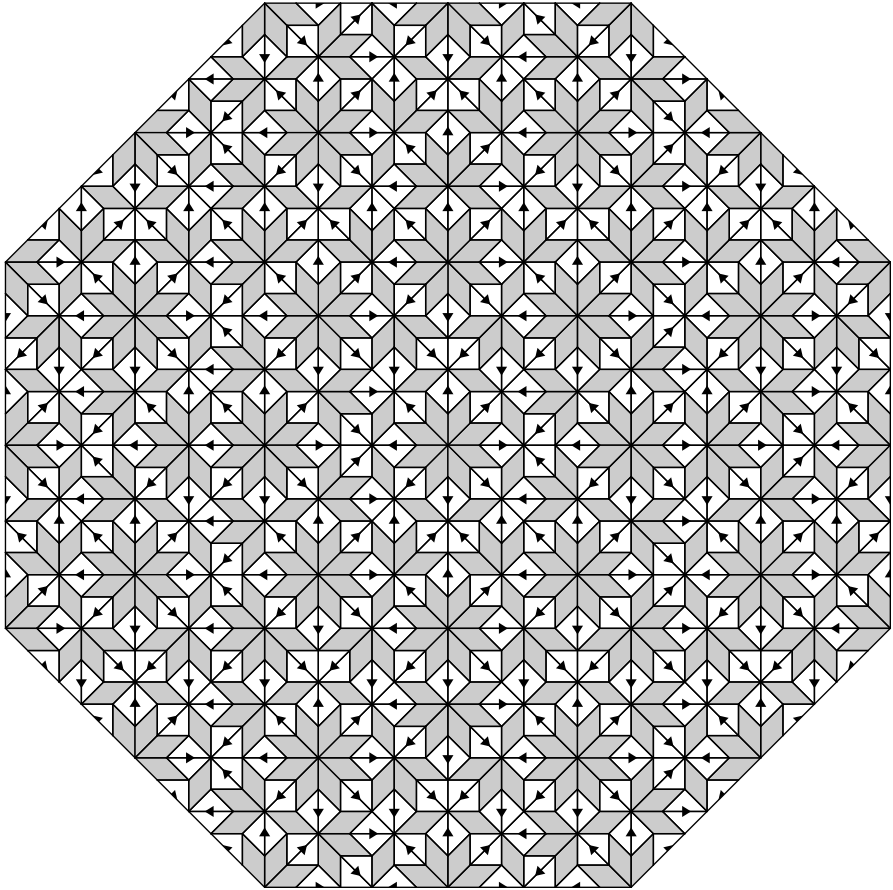}}
\caption{Central patch of an Ammann-Beenker tiling with exact 
$D_{8}$ symmetry.}
\label{fig:absymm}
\end{figure}

In contrast, Figure~\ref{fig:absymm} shows the central patch of a
fixed point tiling (for the same inflation rule, but starting from an
eightfold symmetric seed; the patch shown is obtained by two inflation
steps from the central octagonal patch comprising $16$ triangles and
$16$ rhombuses), this time with complete $D_{8}$ symmetry. Although
these two tilings are distinct, they are locally indistinguishable in
the sense that any finite patch of one occurs in the other (and vice
versa); see below for a precise definition. In fact, any finite patch
re-occurs in a relatively dense fashion. This property is called
\emph{repetitivity}, and is a direct consequence of the underlying
(stone) inflation structure (where stone inflation refers to the
property that the inflated tiles are dissected into complete tiles in
the inflation rule of Figure~\ref{fig:abinfltil}). Both from a
mathematical and from a physical perspective, structures that are
locally indistinguishable should be considered from a unified point of
view; see also the discussion in \cite{Lif}.  In particular, one needs
a suitably adapted notion of symmetry. For example, the diffraction
patterns of our two examples are identical and show perfect $D_{8}$
symmetry; see \cite{BGrev,BG12} for details and an illustration.

\section{Symmetries of LI classes and hulls}

Two Delone sets $\vL$ and $\vL'$ in $\RR^{d}$ are called \emph{locally
  indistinguishable} (or LI for short) if any finite cluster of $\vL$
(which is a non-empty set of the form $\vL\cap K$ for some compact
$K\subset\RR^{d}$) is, up to translation, also a cluster of $\vL'$ and
vice versa. More formally, for any compact $K,K'\subset\RR^{d}$,
there are translations $t,t'\in\RR^{d}$ such that 
\[
   (t+\vL)\cap K \, = \, \vL'\cap K
   \quad \text{and}\quad
   (t'+\vL')\cap K' \, = \, \vL\cap K'\ts .
\]
Note that also the terms locally isomorphic or locally isometric are
in use for this property (but have the same abbreviation). Clearly, LI
is an equivalence relation on the class of Delone sets, and the
equivalence class of a single $\vL$ is called its LI class, denoted by
$\mathrm{LI}(\vL)$.

Beyond LI, one needs a notion to formalise relations between point
sets and tilings, or, more generally, between distinct but related
patterns. A useful concept, called local derivability, was introduced
in \cite{BSJ}; see also \cite{B}. A pattern (where we are thinking
mainly of tilings or point sets here, though more general structures
are possible) $T'$ in $\mathbb{R}^{d}$ is said to be \emph{locally
  derivable} from a pattern $T$ in $\mathbb{R}^{d}$, when a compact
neighbourhood $K\subset\RR^{d}$ of $0$ exists such that, whenever the
$K$-patches (defined as suitable intersections with $K$) of $T$
centred at two positions $x\in\mathbb{R}^{d}$ and $y\in\mathbb{R}^{d}$
coincide, $T'$ also looks the same in the vicinity of $x$ and $y$; we
refer to \cite{BSJ,TAO} for further details and formal aspects.  This
means that there exists a strictly \emph{local} rule to construct $T'$
from $T$.

Next, two patterns $T$ and $T'$ in $\mathbb{R}^{d}$ are called
\emph{mutually locally derivable} or MLD if $T'$ is locally derivable
from $T$ and vice versa. An example of two patterns that are MLD are
the (undecorated) Ammann-Beenker tiling of Figure~\ref{fig:absymm} and
the set of its vertex points --- clearly, knowing the tiling you know
the position of the vertices, while it is slightly less obvious that
knowing the vertices one can locally reconstruct the squares and the
rhombuses, including the arrow direction. Mutual local derivability is
an equivalence concept, and we refer to the corresponding equivalence
classes as MLD classes.  For the Ammann-Beenker tiling of
Figure~\ref{fig:abinfltil} and \ref{fig:absymm}, the set of vertex
points is thus a Delone set representative of the corresponding MLD
class. Note that mutual local derivability consistently extends from
tilings or Delone sets to the corresponding LI classes, wherefore
an MLD class is naturally composed of entire LI classes. 
 \medskip

In our setting, a linear isometry $R\in \mathrm{O}(d,\RR)$ is called
an \emph{LI symmetry} of $\vL\subset\RR^{d}$, if $\vL$ and $R\vL$ are LI,
or (equivalently) if $R\vL\in\mathrm{LI}(\vL)$. Some further thought
reveals that each element of $\mathrm{LI}(\vL)$ shares the same
symmetry in this sense, so that the symmetry is actually a property of
the entire LI class. This definition resolves the problem encountered
with the two Ammann-Beenker tilings (respectively their representing
Delone sets) above, both now having symmetry group $D_{8}$ in this
sense. This is true because the non-symmetric version has the property
that all its rotated copies (by multiples of $\pi/4$) are in the same
LI class. Due to the orientations of the structural elements, no
further symmetries are possible in this example.

\begin{figure}[t]
\centerline{\includegraphics[width=0.8\textwidth]{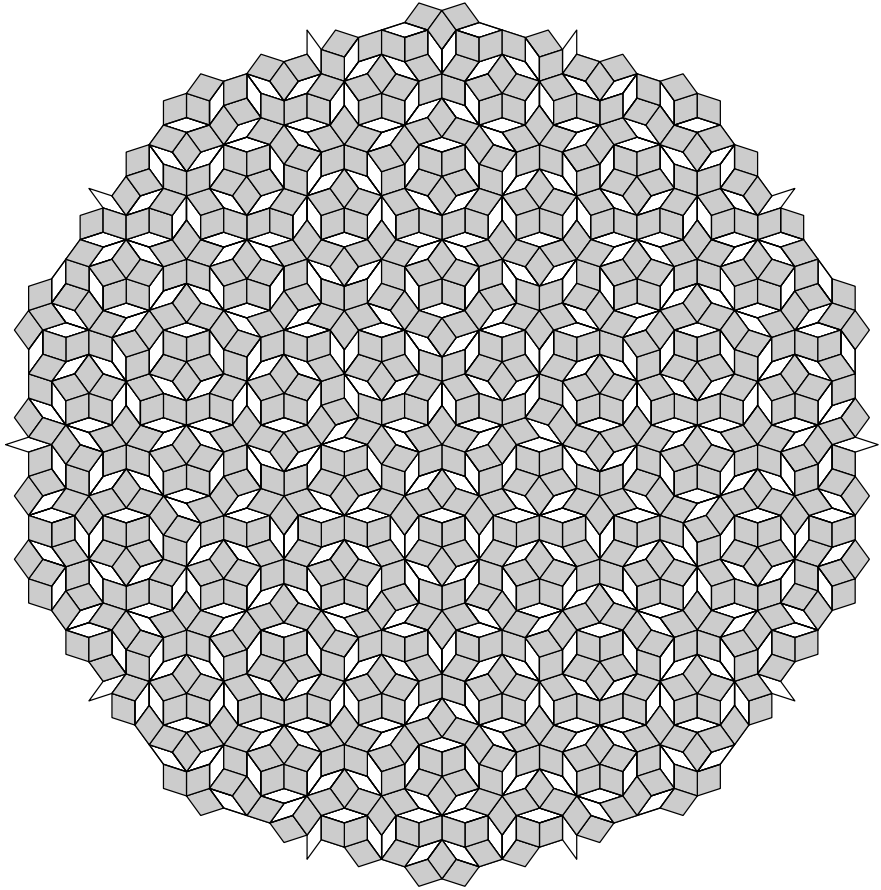}}
\caption{Central patch of a rhombic Penrose tiling with exact $D_{5}$ 
symmetry.}
\label{fig:pen}
\end{figure}

In general, an LI class can have a larger symmetry group than any of
its members individually. For instance, this happens for the rhombic
Penrose tiling. Here, the maximal exact or individual symmetry group
is $D_{5}$, as in the example shown in Figure~\ref{fig:pen}, while any
such Penrose tiling is LI with a copy that is rotated through
$\pi/5$. In particular, the images of any finite patch under the
$D_{10}$ operations occur in the same tiling.  Consequently, the
symmetry group in the LI sense is $D_{10}$, which is also the symmetry
of the diffraction pattern of \emph{any} rhombic Penrose tiling.

This approach is adequate and satisfactory as long as the LI class is
a closed set in the local rubber topology \cite{BL}, where two Delone
sets $\vL$ and $\vL'$ are $\varepsilon$-close when they agree on a
central ball of radius $1/\varepsilon$, possibly after moving the
individual points of one set ($\vL$, say) by at most a distance
$\varepsilon$. This means that each individual point $x\in\vL\cap
B_{1/\varepsilon}(0)$ may be independently displaced to match
$\vL'\cap B_{1/\varepsilon}(0)$, but only within
$B_{\varepsilon}(x)$. This topology is induced by the Hausdorff metric
\cite{RW,Boris}.

In general, however, LI classes are not closed. As a trivial example,
consider the set $\vL=\ZZ\setminus \{0\}$, which may be considered as
a set with a single `defect'. This set is clearly LI with any
translate of itself, and one has
\[
   \mathrm{LI}(\vL) \, = \, \{t+\vL\mid t\in\RR\} \ts .
\]
In other words, the LI class is the orbit of $\vL$ under the
translation action of the group $\RR$. However, this LI class is
\emph{not} closed in the local rubber topology, because
\[
    \lim_{n\to\infty} (\alpha+n+\vL) \, = \, \alpha+\ZZ\ts ,
\]
where $\alpha\in\RR$ is arbitrary but fixed, and $n$ is
integer. Consequently, one finds
\[
   \overline{\mathrm{LI}(\vL)} \, = \,  
    \mathrm{LI}(\vL) \cup \mathrm{LI}(\ZZ)
\]  
with $\mathrm{LI}(\ZZ)= \{\alpha + \ZZ \mid 0\le \alpha<1 \}\simeq
\RR/\ZZ$, which is a closed subset.  Obviously, the different
components behave differently under reflections, although they still
have reflection symmetry in the LI sense. Repeating this exercise with
$\vL=\ZZ\setminus\{0,1,3\}$ leads to an orbit closure with components
of distinct symmetry. For physical aspects such as diffraction, the
defects are immaterial, so that the relevant symmetry is not the
individual symmetry of $\vL$, but rather that of $\ZZ$. The latter is
a closed subset of $\mathrm{LI}(\vL)$. We will come back to this point
below in Section~\ref{sec:prob} in a measure-theoretic setting.

To simplify our further discussion, let us assume that the point set
$\vL\subset\RR^{d}$ has \emph{finite local complexity} (FLC). By this
we mean that, given any compact $K\subset\RR^{d}$, there are only
finitely many distinct clusters $(t+K)\cap\vL$ (with $t\in\RR^{d}$) up
to translations. If $\vL$ is FLC, it must be uniformly discrete. Let
us recall a characterisation of FLC sets due to Schlottmann
\cite{Schl}.

\begin{lemma}
  A point set $\vL\subset\RR^{d}$ is FLC if and only if
  $\vL-\vL=\{x-y\mid x,y\in\vL\}$ is locally\break finite.\hfill$\square$
\end{lemma}

The advantage is that, for FLC sets, we may replace the local rubber
topology by the simpler \emph{local topology} (both agree on FLC sets
\cite{BL}). In the local topology, two FLC sets $\vL$ and $\vL'$ are
$\varepsilon$-close when $\vL'\cap B_{1/\varepsilon}(0) = (t+\vL)\cap
B_{1/\varepsilon}(0)$ for some $t\in B_{\varepsilon}(0)$, so that only
the entire set $\vL$ may be displaced (rather than its points
individually as before).

Given an FLC set $\vL\subset\RR^{d}$, we define its (continuous)
\emph{hull} as 
\[
    \XX(\vL) \, = \, \overline{\{t+\vL\mid t\in\RR^{d}\}}\ts ,
\]
where the closure is taken in the local topology; compare \cite{BHZ}
and references therein. A classic result on $\XX(\vL)$ is the following
\cite{RW,Schl,BHZ}.

\begin{proposition}
  If $\vL\subset\RR^{d}$ is an FLC set, its hull\/ $\XX(\vL)$ is
  compact in the local topology. \hfill$\square$
\end{proposition}

If $\vL\subset\RR^{d}$ is a crystallographic Delone set, with $\vG$ as
its lattice of periods, one obtains the hull 
\[
   \XX(\vL) \, = \, \overline{\{ t+\vL \mid t\in\RR^{d} \}}
  \, = \, \{ t+\vL \mid t\in\mathrm{FD}(\vG) \}
  \, \simeq\, \RR^{d}/\vG \, \simeq\, \mathbb{T}^{d} ,
\]
where $\mathrm{FD}(\vG)$ is a (measurable) fundamental domain of the
lattice $\vG$ and $\mathbb{T}^{d}$ denotes the $d$-torus. This is an
example of a particularly simple hull. In contrast, the hull of the
Ammann-Beenker point set, say, is a rather complicated topological
space, which locally looks like a product of $\RR^{2}$ with a Cantor
set \cite{BHZ,Sadun}.

In many important cases, the hull and the LI class of a Delone set
$\vL$ coincide. This happens if and only if $\vL$ is repetitive
\cite{Schl}.  Examples include the Ammann-Beenker and Penrose vertex
point sets from above. In general, the hull is a larger set than the
LI class. In such a case, one defines the symmetry of $\vL$ via the
hull $\XX(\vL)$, where $R\in \mathrm{O}(d,\RR)$ is a hull symmetry of
$\vL$ if $R\ts\XX(\vL)=\XX(\vL)$. This is the symmetry concept from
dynamical systems theory \cite{Rob,RW}, where the compact space
$\XX(\vL)$ is invariant under the translation action of
$\RR^{d}$. Since translations are used to define the hull, the
built-in translation invariance implies that isometries and
translations have to be treated differently. Let us now turn our
attention to appropriate notions of non-periodicity in this context.

\section{Non-periodicity versus aperiodicity}

Let us define a concept of aperiodicity that is in line with the
notion used in \cite[Ch.~10]{GS}.  As introduced above, a locally
finite point set $\vL\subset\RR^{d}$ is called \emph{non-periodic}
when $\per(\vL)=\{0\}$. The set $\ZZ\setminus\{0\}$ is non-periodic as
a subset of $\RR$.  However, the hull $\XX(\ZZ\setminus\{0\})$
contains $\ZZ$ as an element, which is clearly periodic. This
motivates the following refining definition.

\begin{definition}\label{def1}
  A locally finite point set $\vL\subset\RR^{d}$ is called\/
  \emph{topologically aperiodic}, or\/ \emph{aperiodic} for short,
  when all elements of its hull\/ $\ts\XX(\vL)$ are non-periodic.
\end{definition}

The notion of aperiodicity is thus stronger than that of
non-periodicity. One has the implications
\[
   \text{$\vL\,$ aperiodic} \quad\Rightarrow\quad
   \text{$\vL\,$ non-periodic} \quad\Rightarrow\quad
   \text{$\vL\,$ non-crystallographic}, 
\]
none of which is reversible in general. The Ammann-Beenker and Penrose
point sets from above are aperiodic in this stronger sense. The
importance of the concept of aperiodicity stems from the fact that 
it discards `trivial' cases of non-periodicity, as in our initial
example.\medskip  

\begin{figure}[t]
\centerline{\raisebox{5ex}{\includegraphics[width=0.2\textwidth]{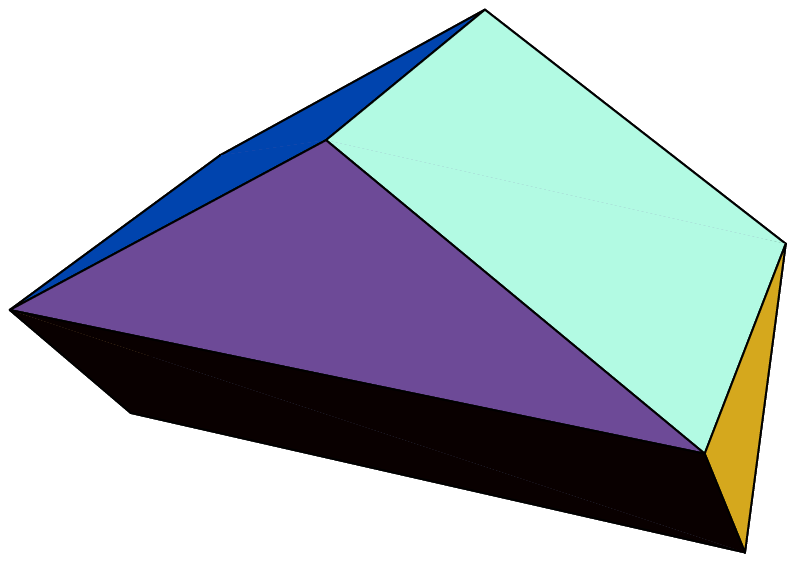}}\hfill
\includegraphics[width=0.8\textwidth]{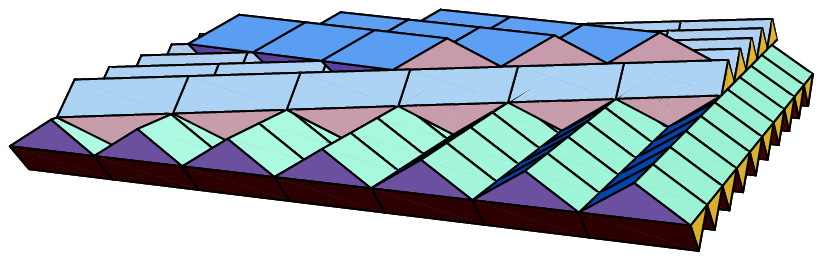}}
\caption{The SCD monotile (left) and the layer structure of the corresponding 
tiling (right).}
\label{scd}
\end{figure}

Let us move on to a three-dimensional structure with interesting
properties, which will lead to a further refinement of our definition
of aperiodicity. Consider a bi-prism with rhombic basis, such as the
one shown in the left panel of Figure~\ref{scd}. Depending on the
angle $\varphi$ between the two `roof ridges', we call the prototile
\emph{commensurate} (for $\varphi\in\pi\mathbb{Q}$) or
\emph{incommensurate} (for $\varphi\not\in\pi\mathbb{Q}$). The
incommensurate SCD prototile was introduced by Schmitt, Conway and
Danzer, see \cite{Danzer,BF} for details and background, as an example
of a monotile or `einstein', which is a single prototile that admits
space-filling tilings, but no tiling with non-trivial periods. Note
that we do not allow the use of the reflected tile here. As before,
one can replace any tiling by a suitable Delone set (for instance via
the vertices of the Delone complex, which is the dual of the Voronoi
complex).  For convenience, we use the tiling language for this
example, but phrase results for a suitable Delone representative of
the corresponding MLD class.

Joining SCD tiles face to face at the triangular facets, one can
(only) form a closed planar layer, which is periodic with a
two-dimensional lattice of periods $\vG$. The top of this layer shows
ridges and valleys, all parallel to each other, while the bottom shows
the corresponding structure `upside down', with the directions of the
ridges related by an angle $\varphi$. In order to continue the tiling
process, one has to stack such layers as indicated in the right panel
of Figure~\ref{scd}. Each layer is rotated with respect to the
previous layer by the angle $\varphi$. For a finite number of layers,
it is still possible that non-trivial periods persist, which happens
precisely if the rotation $R_{\varphi}$ is a coincidence rotation of
the planar lattice of periods of a single layer, so when $\vG$ and
$R_{\varphi}\vG$ share a common sublattice of finite index. However,
with each layer added, the common sublattice gets sparser, and
$\bigcap_{m\in\ZZ}R^{m}_{\varphi}\vG = \{0\}$. Consequently, a space-filling
tiling built from an incommensurate SCD tile cannot have any
non-trivial periods, and hence is non-periodic.

In general, an SCD tiling is not repetitive, even when this property
is considered with respect to Euclidean motions, which is appropriate
here. Repetitive SCD tilings can be obtained \cite{BF} if
$\cos(\varphi)\in\mathbb{Q}$. For a repetitive tiling, the
corresponding hull is again a compact set, and contains only congruent
copies of a single SCD tiling and limits of convergent sequences of
congruent copies. This construction leaves no freedom to generate
non-trivial periods, which implies the following result according to
Definition~\ref{def1}.

\begin{proposition}\label{prop3}
  Consider a repetitive, incommensurate SCD tiling, and let $\vL$ be a
  Delone representative of its MLD class. Then, $\vL$ is aperiodic.
  \hfill$\square$
\end{proposition}

Clearly, one then also calls the corresponding tiling itself
aperiodic, after the obvious extension of our above definitions to
tilings and (more generally) to patterns; see \cite{TAO} for details.
Consequently, the rhombic Penrose tiling is aperiodic, as is the
Ammann-Beenker tiling (and many others). 
\medskip

\begin{figure}[t]
\centerline{\includegraphics[width=0.5\textwidth]{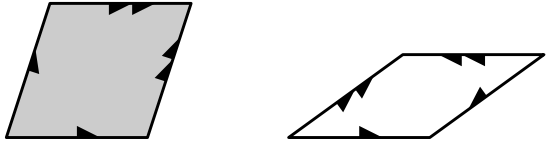}}
\caption{Aperiodic prototile set for the rhombic Penrose tiling.}
\label{penmatch}
\end{figure}

The SCD tile is an example of an \emph{aperiodic prototile set}, a
term frequently met in the literature; compare \cite{GS}. This refers
to a (usually finite) set of prototiles (possibly with markers) which
admit space-filling tilings, but only non-periodic ones. An example of
an aperiodic prototile set are the two decorated Penrose rhombuses of
Figure~\ref{penmatch}, with the imposed local rules (or matching
rules) that arrows on adjacent edges have to match (and the tiling has
to be face to face). Note that the two `naked' Penrose rhombuses
(without arrow decorations) fail to form an aperiodic set (because
they are obviously compatible with periodic tilings), while the
decorated versions of Figure~\ref{penmatch} do the job, via these
arrow matching rules \cite{GS}, in the sense that any space-filling
tiling of the plane that obeys the arrow matching rules is MLD with
the Penrose tiling of Figure~\ref{fig:pen}.  Note that in any infinite
rhombic Penrose tiling, the arrow decorations of the edges can be
reconstructed from the undecorated version in a local fashion, so the
decorated and undecorated tilings are indeed MLD.

This property is not met by most of the other tilings of the aperiodic
tiling zoo, including the Ammann-Beenker tiling. The latter also
possesses aperiodic local rules, which can be realised as additional
decorations, but these cannot be locally reconstructed from the
undecorated tiling (while removing the additional decorations clearly
provides a local rule in the converse direction). In that case, the
decorated and undecorated tilings form different MLD classes; see
\cite{Franz} for a proof of this statement.

Given an aperiodic prototile set, all admissible space-filling tilings
together form a space $\XX$, wherefore our above notion of
aperiodicity fits this situation as well. However, note that a tiling
can be aperiodic in our setting without emerging from an aperiodic
prototile set, and the undecorated Ammann-Beenker tiling is an example
(as are any one-dimensional aperiodic tilings, since aperiodicity
cannot be enforced by local rules in one dimension). Note that, for
the Ammann-Beenker tiling, an aperiodic prototile set with additional
markings exist \cite{GS}, yielding the decorated Ammann-Beenker
tilings from which the undecorated Ammann-Beenker tilings can be
locally derived (by removing the additional markings), but not vice
versa, as the additional markings cannot be locally recovered. A
similar result holds more generally for inflation tilings
\cite{Chaim}, in the sense that there is an, in principle,
constructive method to add local information to obtain an aperiodic
prototile set; however, this usually leads to a much larger number of
prototiles and is thus not used in practice.

As this discussion shows, the meaning of `aperiodicity' depends on the
context, but our main point here is that the notion built on the LI
class (rather than on an individual tiling or Delone set) is
consistent.  \medskip

Returning to the SCD tiling, the result of Proposition~\ref{prop3} is
perhaps not entirely satisfactory, since an individual SCD tiling of
this type may still possess a symmetry in the form of a screw
axis. Even though the tiling admits no non-trivial translation
symmetry, its symmetry group thus contains an infinite subgroup. The
following definition introduces a stronger notion that separates such
structures from examples like the Ammann-Beenker or Penrose tilings
considered above, where the remaining symmetry groups are always
finite.

\begin{definition}
  A locally finite point set $\vL\subset\RR^{d}$ is called\/
  \emph{strongly aperiodic} when it is aperiodic and when the
  individual symmetry group of each element of\/ $\XX(\vL)$ is a finite
  group.
\end{definition}

\begin{figure}[t]
\centerline{\includegraphics[width=0.9\textwidth]{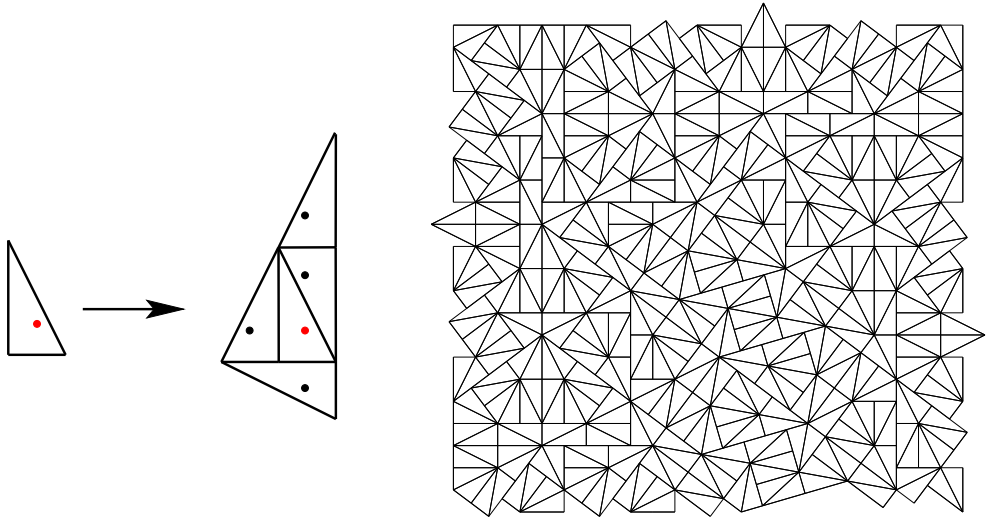}}
\caption{Inflation rule (left) and patch (right) of the pinwheel
  tiling. The prototiles are right triangles with side length $1$, $2$
  and $\sqrt{5}$.  The dots represent control points of an equivalent
  Delone set. Starting from a single triangle with the control point
  in the origin leads to a fixed point tiling with circular symmetry.}
\label{fig:pin}
\end{figure}

With this notion, repetitive incommensurate SCD tilings are aperiodic,
but not strongly aperiodic. Note that it is possible to have strong
aperiodicity in the presence of continuous symmetries of the hull. An
example is provided by the pinwheel tiling \cite{Rad}. This is a
tiling of the plane with a single triangular prototile (and its
reflected version), defined by the inflation rule shown in the left
panel of Figure~\ref{fig:pin}, where the dots indicate a possible
choice of representative (or control) points for an MLD Delone set
\cite{BFG}. The reflected tile is inflated using the reflected rule.
The inflation rule introduces an angle $2\arctan(\frac{1}{2})$ which
is incommensurate with $\pi$, so a 'new' direction is introduced in
each inflation step, so a fixed point tiling shows complete circular
symmetry (while and individual tiling has at most a reflection
symmetry), and as a consequence has circularly symmetric diffraction
\cite{MPS,BFG} (although a complete description of the diffraction
measure is still lacking).  A small patch of a pinwheel tiling is
shown in the right panel of Figure~\ref{fig:pin}. Various
generalisation of pinwheel tilings exist \cite{Fre08}, including
higher-dimensional `quaquaversal' tilings \cite{CR98}; see also
\cite[Ch.~4.3]{RadBook}.

\section{Probabilistic extensions and outlook}\label{sec:prob}

The notion of a hull is an example of an \emph{ensemble} of point
sets, which is a well-known concept in statistical physics, too.  For
instance, consider the (discrete) ensemble $\XX$ of all subsets of
$\ZZ^{d}$ (which can later be extended to a continuous ensemble by
translations). This ensemble can be identified with the compact set
$\{0,1\}^{\ZZ^{d}}$ and may be considered as a model for a lattice
gas. In ergodic theory, this is known as the full shift space
\cite{Walters}. Note that $\XX$ is the orbit closure (under the action
of $\ZZ^{d}$) of a typical element of $\XX$, but \emph{not} of every
element. More precisely, there exists a probability measure $\mu$ (in
fact, more than one) on $\XX$ that quantifies this statement in a
certain way. One possibility is conveniently described via the
`obvious' Bernoulli measure on $\XX$. For $0<p<1$, the latter is
defined via the measure of cylinder sets, giving probability $p$ for
$1$ (occupied) and $1-p$ for $0$ (empty) per site; see \cite{Walters}
for details on the construction. Any such measure is
$\ZZ^{d}$-invariant by construction. Consequently, one obtains a
measure-theoretic dynamical system $(\XX,\ZZ^{d},\mu)$.

In such a dynamical system, $\mu$-almost all elements are
non-periodic, while $\XX$ also contains lots of periodic or even
crystallographic elements (such as the completely occupied
lattice). However, all periodic elements together are still a subset
of $\XX$ of measure $0$. This motivates the following probabilistic
counterpart of Definition~\ref{def1}.

\begin{definition}
  A measure-theoretic dynamical system\/ $(\XX,\RR^{d},\mu)$ is
  called\/ \emph{statistically (or metrically) aperiodic} when the
  subset of elements of\/ $\XX$ with non-trivial periods is a null set
  for the measure\/ $\mu$.
\end{definition}

The term `metric' is often used in number theory or dynamical systems
theory to refer to a situation where certain properties hold almost
surely relative to a reference measure.

Returning to our example $\vL=\ZZ\setminus\{0\}$ from above, one can
show that the corresponding hull $\XX(\vL)$ admits only one invariant
probability measure.  This has support on $\mathrm{LI}(\ZZ)$ and hence
gives no weight at all to the subset $\mathrm{LI}(\ZZ\setminus\{0\})$.
This means that, in the probabilistic setting, the dynamical systems
point of view tells us that $\vL$ should be considered as an
essentially periodic system, and not as a non-periodic one --- quite
in agreement with the physical intuition about the role of a single
defect in an otherwise periodic structure.

It is clear that one can use this probabilistic setting to also extend
the notion of symmetry, then often referred to as \emph{statistical
  symmetry}. There are several possible variants of it. For any given
Delone set with well-defined patch frequencies, one could speak of a
statistical rotation or reflection symmetry when the image of each
patch of the set occurs with the same frequency. For instance, each
Penrose point set from above has statistical $D_{10}$ symmetry in this
sense. This coincides with the symmetry of the LI class, because the
latter is closed and uniquely ergodic as a dynamical system \cite{Rob}
under the action of $\RR^{2}$.

A larger class would comprise equilibrium systems, which essentially
are ensembles together with a translation invariant Gibbs measure; see
\cite{Geo} for a systematic exposition. Such systems, and
probabilistic methods in general, are presently receiving growing
attention, both theoretically and practically; see \cite{BG12} and
references therein. It is thus desirable and necessary to further
extend the classic notions of crystallography to this viewpoint. We
believe that this will be vital for a better understanding of
important stochastic system such as random tilings or disordered
structures.

\section*{Acknowledgement}

This work was supported by the German Research Council
(DFG), within the CRC 701.

\begin{small}

\end{small}

\end{document}